\documentclass{amsart}

\usepackage{palatino}
\usepackage{amssymb}
\usepackage{amsmath}
\usepackage{amsthm}
\usepackage[english]{babel}

\newcommand{\N}{\mathbb{N}}

\newtheorem{theorem}{Theorem}[section]
\newtheorem{proposition}[theorem]{Proposition}
\newtheorem{lemma}[theorem]{Lemma}
\newtheorem{corollary}[theorem]{Corollary}
\newtheorem{definition}[theorem]{Definition}
\newtheorem{example}[theorem]{Example}
\newtheorem{remark}[theorem]{Remark}

\title{Representation of numerical semigroups by Dyck paths}

\author{Maria Bras-Amor\'os and Anna de Mier}
\address[Maria Bras-Amor\'os]{Departament d'Enginyeria de la
  Informaci\'o i de les Comunicacions, 
Universitat Aut\`onoma de Barcelona, 
08193 Bellaterra, Spain} \email { maria.bras@uab.cat
} 
\address[Anna de Mier]{Departament de Llenguatges i Sistemes
Inform\`atics, Universitat Polit\`ecnica de Catalunya, Jordi
Girona 1--3, 08034 Barcelona, Spain} \email { anna.de.mier@upc.edu
} 

\begin{document}

\maketitle

\begin{abstract}
We introduce square diagrams that represent numerical semigroups
and we obtain an injection from the set of numerical semigroups
into the set of Dyck paths.
\end{abstract}

\section*{Introduction}
A {\it numerical semigroup} is a subset of $\N_0$ closed under
addition and with finite complement in $\N_0$
\cite{RoGa,HoLiPe:agc}. The {\it genus} of a numerical semigroup
$\Lambda$ is the number of elements in $\N_0\setminus \Lambda$,
which are called \emph{gaps}. A {\it Dyck path} of order $n$ is a
lattice path from $(0,0)$ to $(n,n)$ consisting of up-steps
$\uparrow=(0,1)$ and right-steps $\rightarrow=(1,0)$ and never
going below the diagonal $x=y$.

We introduce the notion of the {\it square diagram} of a numerical
semigroup and analyze some properties of numerical semigroups such as
their weight or symmetry by means of the square diagram.

We prove that any numerical semigroup is represented by a unique
Dyck path of order given by its genus.

\section{Square diagram}
Given a numerical semigroup $\Lambda$ define $\tau(\Lambda)$ as
the path with origin $(0,0)$ and steps $e(i)$ given by
$$e(i)=
\left\{
\begin{array}{ll}
\rightarrow & \mbox{ if } i\in \Lambda,\\
\uparrow & \mbox{ if } i\not\in\Lambda,\\
\end{array}
\right. \qquad \mbox{ for } 1\leq i\leq 2g.$$ We denote it as the
{\it square diagram} of $\Lambda$. For instance, the set
$$\{0, 4, 8, 12, 16, 17, 19, 20, 21, 23, 24, 25, 27, 28, 29,
31\}\cup\{i\in\N_0: i\geq 31\}$$
is a numerical semigroup and its square diagram is the following one:

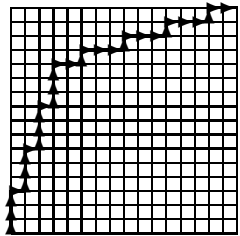
\setlength{\unitlength}{0.187500cm}
\begin{center}
\begin{picture}(16.,16.)

\put(0,0){\line(1,0){16.}}
\put(0,0){\line(0,1){16.}}
\put(0,1){\line(1,0){16.}}
\put(1,0){\line(0,1){16.}}
\put(0,2){\line(1,0){16.}}
\put(2,0){\line(0,1){16.}}
\put(0,3){\line(1,0){16.}}
\put(3,0){\line(0,1){16.}}
\put(0,4){\line(1,0){16.}}
\put(4,0){\line(0,1){16.}}
\put(0,5){\line(1,0){16.}}
\put(5,0){\line(0,1){16.}}
\put(0,6){\line(1,0){16.}}
\put(6,0){\line(0,1){16.}}
\put(0,7){\line(1,0){16.}}
\put(7,0){\line(0,1){16.}}
\put(0,8){\line(1,0){16.}}
\put(8,0){\line(0,1){16.}}
\put(0,9){\line(1,0){16.}}
\put(9,0){\line(0,1){16.}}
\put(0,10){\line(1,0){16.}}
\put(10,0){\line(0,1){16.}}
\put(0,11){\line(1,0){16.}}
\put(11,0){\line(0,1){16.}}
\put(0,12){\line(1,0){16.}}
\put(12,0){\line(0,1){16.}}
\put(0,13){\line(1,0){16.}}
\put(13,0){\line(0,1){16.}}
\put(0,14){\line(1,0){16.}}
\put(14,0){\line(0,1){16.}}
\put(0,15){\line(1,0){16.}}
\put(15,0){\line(0,1){16.}}
\put(0,16){\line(1,0){16.}}
\put(16,0){\line(0,1){16.}}

{\thicklines
\put(0,0){{\vector(0,1){1.}}}
\put(0,1){{\vector(0,1){1.}}}
\put(0,2){{\vector(0,1){1.}}}
\put(0,3){{\vector(1,0){1.}}}
\put(1,3){{\vector(0,1){1.}}}
\put(1,4){{\vector(0,1){1.}}}
\put(1,5){{\vector(0,1){1.}}}
\put(1,6){{\vector(1,0){1.}}}
\put(2,6){{\vector(0,1){1.}}}
\put(2,7){{\vector(0,1){1.}}}
\put(2,8){{\vector(0,1){1.}}}
\put(2,9){{\vector(1,0){1.}}}
\put(3,9){{\vector(0,1){1.}}}
\put(3,10){{\vector(0,1){1.}}}
\put(3,11){{\vector(0,1){1.}}}
\put(3,12){{\vector(1,0){1.}}}
\put(4,12){{\vector(1,0){1.}}}
\put(5,12){{\vector(0,1){1.}}}
\put(5,13){{\vector(1,0){1.}}}
\put(6,13){{\vector(1,0){1.}}}
\put(7,13){{\vector(1,0){1.}}}
\put(8,13){{\vector(0,1){1.}}}
\put(8,14){{\vector(1,0){1.}}}
\put(9,14){{\vector(1,0){1.}}}
\put(10,14){{\vector(1,0){1.}}}
\put(11,14){{\vector(0,1){1.}}}
\put(11,15){{\vector(1,0){1.}}}
\put(12,15){{\vector(1,0){1.}}}
\put(13,15){{\vector(1,0){1.}}}
\put(14,15){{\vector(0,1){1.}}}
\put(14,16){{\vector(1,0){1.}}}
\put(15,16){{\vector(1,0){1.}}}
}

\end{picture}\end{center}

The {\it conductor} of a numerical semigroup is the unique element
$c$ in $\Lambda$ such that $c-1\not\in\Lambda$ and
$c+\N_0\subseteq \Lambda$. The \emph{enumeration} of $\Lambda$ is
the unique bijective increasing map
$\lambda:\N_0\longrightarrow\Lambda$. We use $\lambda_i$ to denote
$\lambda(i)$. The $ith$ {\it partial genus} may be defined as
\begin{center}
$g(i)=\lambda_i-i=\#$
gaps smaller than $\lambda_i$.
\end{center}

Note that
the following statements are satisfied:
\begin{enumerate}
\item $g(0)=0$,
\item $g(i)\leq g(i+1)$,
\item $g(i)=g$ for all $i\geq \lambda^{-1}(c)$,
\item $g(i)=g$ for all $i\geq g$ (consequence of (iii)).
\end{enumerate}

The poitns with integer coordinates in the square diagram of
$\Lambda$ are all points in
$$\{(i,g(i)): 0\leq i\leq g\}\cup
\{(i-1,g(i)): 1\leq i\leq g\}$$
together with
the points contained in the vertical lines
from $(i-1,g(i-1))$
to $(i-1,g(i))$ whenever $g(i-1)<g(i)$.
In particular,
$\tau(\Lambda)$ goes from $(0,0)$ to $(g,g)$.
So it is included in the square grid from $(0,0)$ to
$(g,g)$.
This is why we call this diagram the square diagram of $\Lambda$.

\section{Square diagram of symmetric semigroups}

The conductor $c$ of any numerical semigroup
satisfies $c\leq 2g$, where $g$ is the genus of the semigroup.
When $c=2g$ the numerical semigroup is said to be {\it symmetric}
\cite{KiPe:telescopic,HoLiPe:agc}.
It is well known that all semigroups
generated by two integers
are symmetric. As a consequence of
the definition of symmetric semigroups we have the following proposition:

\begin{proposition}A numerical semigroup $\Lambda$
is symmetric if and only if its square diagram satisfies
$e(2g-1)=\uparrow$.
\end{proposition}

The next proposition is a well known result on symmetric
semigroups:

\begin{proposition}
A numerical semigroup $\Lambda$
with conductor $c$
is symmetric if and only if
for any non-negative integer $i$, if $i$ is a gap,
then $c-1-i$ is a non-gap.
\end{proposition}

The proof
can be found in
\cite[Remark~4.2]{KiPe:telescopic}
and \cite[Proposition~5.7]{HoLiPe:agc}.
It follows by counting the number of gaps and non-gaps
smaller than the conductor and the fact that
if $i$ is a non-gap then $c-1-i$ must be a
gap because otherwise $c-1$ would also be a non-gap.
As a consequence we have the following property on the square diagram
of symmetric semigroups:

\begin{corollary}
A numerical semigroup with genus $g$ is symmetric if and only if
the intersection of its square diagram with the square determined
by the points $(0,0)$ and $(g-1,g-1)$ is symmetric with respect to
the diagonal from $(0,g-1)$ to $(g-1,0)$.
\end{corollary}

\begin{example}
The set
$$\{0, 4, 8, 12, 16, 17, 18, 20, 21, 22, 24, 25, 26, 28, 29, 30, 32\}
\cup\{i\in\N_0:i\geq 32\}$$
is a numerical
 semigroup and it is symmetric.
Its square diagram is therefore symmetric
with respect to
the diagonal from
$(0,g-1)$ to $(g-1,0)$.

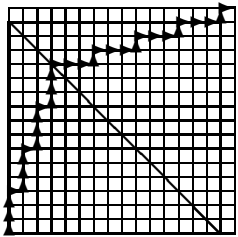
\setlength{\unitlength}{0.187500cm}
\begin{center}
\begin{picture}(16.,16.)

\put(0,0){\line(1,0){16.}}
\put(0,0){\line(0,1){16.}}
\put(0,1){\line(1,0){16.}}
\put(1,0){\line(0,1){16.}}
\put(0,2){\line(1,0){16.}}
\put(2,0){\line(0,1){16.}}
\put(0,3){\line(1,0){16.}}
\put(3,0){\line(0,1){16.}}
\put(0,4){\line(1,0){16.}}
\put(4,0){\line(0,1){16.}}
\put(0,5){\line(1,0){16.}}
\put(5,0){\line(0,1){16.}}
\put(0,6){\line(1,0){16.}}
\put(6,0){\line(0,1){16.}}
\put(0,7){\line(1,0){16.}}
\put(7,0){\line(0,1){16.}}
\put(0,8){\line(1,0){16.}}
\put(8,0){\line(0,1){16.}}
\put(0,9){\line(1,0){16.}}
\put(9,0){\line(0,1){16.}}
\put(0,10){\line(1,0){16.}}
\put(10,0){\line(0,1){16.}}
\put(0,11){\line(1,0){16.}}
\put(11,0){\line(0,1){16.}}
\put(0,12){\line(1,0){16.}}
\put(12,0){\line(0,1){16.}}
\put(0,13){\line(1,0){16.}}
\put(13,0){\line(0,1){16.}}
\put(0,14){\line(1,0){16.}}
\put(14,0){\line(0,1){16.}}
\put(0,15){\line(1,0){16.}}
\put(15,0){\line(0,1){16.}}
\put(0,16){\line(1,0){16.}}
\put(16,0){\line(0,1){16.}}

{\thicklines
\put(0,0){{\line(1,0){15.}}}
\put(0,0){{\line(0,1){15.}}}
\put(0,15){{\line(1,0){15.}}}
\put(15,0){{\line(0,1){15.}}}
\put(0,15){{\line(1,-1){15.}}}
}

{\thicklines
\put(0,0){{\vector(0,1){1.}}}
\put(0,1){{\vector(0,1){1.}}}
\put(0,2){{\vector(0,1){1.}}}
\put(0,3){{\vector(1,0){1.}}}
\put(1,3){{\vector(0,1){1.}}}
\put(1,4){{\vector(0,1){1.}}}
\put(1,5){{\vector(0,1){1.}}}
\put(1,6){{\vector(1,0){1.}}}
\put(2,6){{\vector(0,1){1.}}}
\put(2,7){{\vector(0,1){1.}}}
\put(2,8){{\vector(0,1){1.}}}
\put(2,9){{\vector(1,0){1.}}}
\put(3,9){{\vector(0,1){1.}}}
\put(3,10){{\vector(0,1){1.}}}
\put(3,11){{\vector(0,1){1.}}}
\put(3,12){{\vector(1,0){1.}}}
\put(4,12){{\vector(1,0){1.}}}
\put(5,12){{\vector(1,0){1.}}}
\put(6,12){{\vector(0,1){1.}}}
\put(6,13){{\vector(1,0){1.}}}
\put(7,13){{\vector(1,0){1.}}}
\put(8,13){{\vector(1,0){1.}}}
\put(9,13){{\vector(0,1){1.}}}
\put(9,14){{\vector(1,0){1.}}}
\put(10,14){{\vector(1,0){1.}}}
\put(11,14){{\vector(1,0){1.}}}
\put(12,14){{\vector(0,1){1.}}}
\put(12,15){{\vector(1,0){1.}}}
\put(13,15){{\vector(1,0){1.}}}
\put(14,15){{\vector(1,0){1.}}}
\put(15,15){{\vector(0,1){1.}}}
\put(15,16){{\vector(1,0){1.}}}
}

\end{picture}
\end{center}

\end{example}

\begin{remark}
There exist paths from $(0,0)$ to $(g-1,g-1)$
which are symmetric with respect to
the diagonal from
$(0,g-1)$ to $(g-1,0)$ but
which do not correspond to a numerical semigroup.
For example,
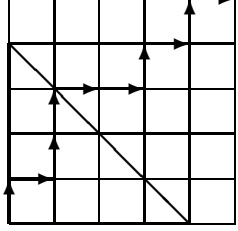
\begin{center}
\setlength{\unitlength}{0.6cm}
\begin{picture}(5.,5.)
\put(0,0){\line(0,1){5}}
\put(1,0){\line(0,1){5}}
\put(2,0){\line(0,1){5}}
\put(3,0){\line(0,1){5}}
\put(4,0){\line(0,1){5}}
\put(5,0){\line(0,1){5}}
\put(0,0){\line(1,0){5}}
\put(0,1){\line(1,0){5}}
\put(0,2){\line(1,0){5}}
\put(0,3){\line(1,0){5}}
\put(0,4){\line(1,0){5}}
\put(0,5){\line(1,0){5}}
{\thicklines
\put(0,0){{\vector(0,1){1}}}
\put(0,1){{\vector(1,0){1}}}
\put(1,1){{\vector(0,1){1}}}
\put(1,2){{\vector(0,1){1}}}
\put(1,3){{\vector(1,0){1}}}
\put(2,3){{\vector(1,0){1}}}
\put(3,3){{\vector(0,1){1}}}
\put(3,4){{\vector(1,0){1}}}
\put(4,4){{\vector(0,1){1}}}
\put(4,5){{\vector(1,0){1}}}
}
{\thicklines
\put(0,0){{\line(0,1){4}}}
\put(0,0){{\line(1,0){4}}}
\put(4,0){{\line(0,1){4}}}
\put(0,4){{\line(1,0){4}}}
\put(0,4){{\line(1,-1){4}}}
}
\end{picture}
\end{center}
This diagram does not correspond to a numerical semigroup
because otherwise, $\lambda_1$ would be $2$, but $4$ would
not belong to the semigroup.
\end{remark}

\section{Weight of a semigroup}

The notion of the weight of a numerical semigroup
has been widely used in the context of Weierstrass semigroups
\cite{Oliveira,CaTo}.

\begin{definition}
Let $\Lambda$ be a numerical semigroup
with genus $g$
and let $l_1,\dots,l_g$ be its gaps.
The {\it weight} of $\Lambda$
is the sum
$$\sum_{i=1}^g \left(l_i-i\right).$$
\end{definition}

In a sense, the weight measures how complicated the semigroup is.
For example, the simplest semigroup is that with gaps $1,2,\dots,g$
and it has weight~$0$.

\begin{proposition}
The weight of a numerical semigroup
is equal to the area over the path in the
square diagram of the semigroup.
\end{proposition}

\begin{proof}
Let $\lambda$ be the enumeration of the semigroup.
The area below the path is equal to the sum
$\sum_{i=1}^g g(i)$,
while the total area of the square diagram is $g^2$.
So, it is enough to prove that
$$\sum_{i=1}^g g(i)+
  \sum_{i=1}^g \left(l_i-i\right)=g^2.$$
But
$\sum_{i=1}^g g(i)+
  \sum_{i=1}^g \left(l_i-i\right)=
\sum_{i=1}^g \left(\lambda_i-i\right)+
  \sum_{i=1}^g \left(l_i-i\right)=
\sum_{i=1}^{2g}i-2\sum_{i=1}^g i=
\sum_{i=g+1}^{2g}i-\sum_{i=1}^g i=
\sum_{i=1}^{g}\left(g+i\right)-\sum_{i=1}^g i=g^2$.

\end{proof}

\section{The square diagram of a numerical semigroup represents a Dyck path}

\begin{lemma}
\label{lemma1}
Let $\Lambda$ be a numerical semigroup with genus $g$ and enumeration $\lambda$.
If $g(i)<i$ then $\lambda_{i+1}=\lambda_i+1$.
\end{lemma}

\begin{proof}
If $g(i)<i$ then there are more non-gaps than gaps in the interval
$[1,\lambda_i]$. So by the Pigeonhole Principle there must be at
least one pair $a,b\in\Lambda$ with $a+b=\lambda_i+1$. Therefore,
$\lambda_{i+1}=\lambda_i+1$.
\end{proof}

\begin{lemma}
\label{lemma2}
Let $\Lambda$ be a numerical semigroup with genus $g$, conductor $c$ and enumeration $\lambda$.
If $g(i)<i$ then $\lambda_i\geq c$.
\end{lemma}

\begin{proof}
Let us show by induction that $\lambda_i+k\in\Lambda$ for all $k\geq 0$.
It is obvious for $k=0$. If $\lambda_i+k'\in\Lambda$ for all $0\leq k'\leq k$ then
$g(i+k)=g(i)<i\leq i+k$ and, by Lemma~\ref{lemma1}, $\lambda_i+(k+1)\in\Lambda$.
\end{proof}

\begin{theorem}
The path $\tau(\Lambda)$
associated to a
numerical semigroup $\Lambda$
is a Dyck path.
\end{theorem}

\begin{proof}
Let $g$ and $c$ be the genus and the conductor of $\Lambda$. It is
enough to show that for all $i$ with $0\leq i \leq g$ we have
$g(i)\geq i$. Indeed, if $g(i)<i$, by Lemma~\ref{lemma2},
$\lambda_i\geq c$ and $g(i)=g$, so $i>g$, a contradiction.
%Since $g(i)=\lambda_i-i$, this is equivalent to showing that
%$\lambda_i\geq 2i$ for all $i$ with $0\leq i\leq g$. Since for $g=0$ there is
%nothing to prove, assume $g\geq 1$. The statement clearly holds for $i=0,1$.
%Assume the result does not hold and let $j+1\leq g$ be the first place where
% the path
%crosses the diagonal; that is, suppose that $\lambda_{j+1}=2j+1$. We show that
%in this situation $2j+k\in \Lambda$ for all $k\geq 0$, hence proving that
%the genus $g$ is $j$, a contradiction. We prove this by induction on $k$.
%The cases $k=0$ and $k=1$  hold by the choice of $j$.
%Assume that $2j+k'$ belongs to $\Lambda$ for all $k'<k$. There are
%$\lfloor {2j+k}/2 \rfloor$ pairs $(a,b)$ with $a\leq b\leq 2j+k-1$ such that
%$a+b=2j+k$  ($a$ and $b$ could be equal); also, only $j$ of the integers
%in the interval $[1,2j+k-1]$
%do not belong to $\Lambda$. Since $\lfloor {2j+k}/2 \rfloor \geq j+1$, at least
%for one of the pairs $(a,b)$ as above both $a$ and $b$ belong to $\Lambda$
%(pigeonhole principle), and
%hence $2j+k=a+b\in \Lambda$.
\end{proof}

\begin{corollary}\label{cor:fitacatalan}
There are at most $\frac{1}{g+1}\binom{2g}{g}$  numerical semigroups of genus $g$.

There are at most $\binom{g-1}{\lceil {g-1}/2 \rceil}$
symmetric numerical semigroups of genus $g$.
\end{corollary}

\begin{proof}
The first assertion is a consequence
of the fact that the number of Dyck paths of order $n$ is given by the Catalan number
$C_n=\frac{1}{n+1}\binom{2n}{n}$
(see~\cite[Corollary 6.2.3(v)]{ec2}).

For the second, recall that the square diagram of a symmetric
numerical semigroup is  symmetric in the square determined by
$(0,0)$ and $(g-1,g-1)$.  Such a symmetric path is determined by
the first half of the path. The result follows from the fact
(see~\cite{cat}) that the number of paths with $m$ steps that
start at $(0,0)$ and never go below the line $x=y$ is
$\binom{m}{\lceil {m}/2 \rceil}$.
\end{proof}

%\bibliographystyle{plain}
%\bibliography{bib}

\end{document}